# Dyck Paths in Four-Dimensional Space

Gennady Eremin
ergenns@gmail.com

January 14, 2019

**Abstract.** In analyzing balanced parentheses, we consider a group of related variables in Dyck paths. In the four-dimensional space, the Dyck triangle is constructed – an integer lattice with Dyck paths.

## 1 Introduction

In multidimensional spaces (dimension 4 and more), virtual constructions are usually described, not real objects. Considering balanced parentheses, we found a close relationship of known combinatorial objects. The Dyck triangle with diagonal paths [1, §2.5], the Dyck triangle with monotonic paths in a square grid [2], and the convolution of the Catalan's matrix [3] are different projections of the same four-dimensional body from the first (positive) hyperoctant.

The above objects are enumerated by Catalan numbers (A000108), these numbers are manifested in various problems of discrete mathematics [4]. We have borrowed Figure 1 from Tom Davis [2]. In the left picture we see a diagonal path (or a mountain range) of 12 links located above the zero line. The right picture shows the corresponding monotonic path in a square grid 6 by 6. The monotonic path leads from the lower left corner to the upper right corner, without crossing the main diagonal.

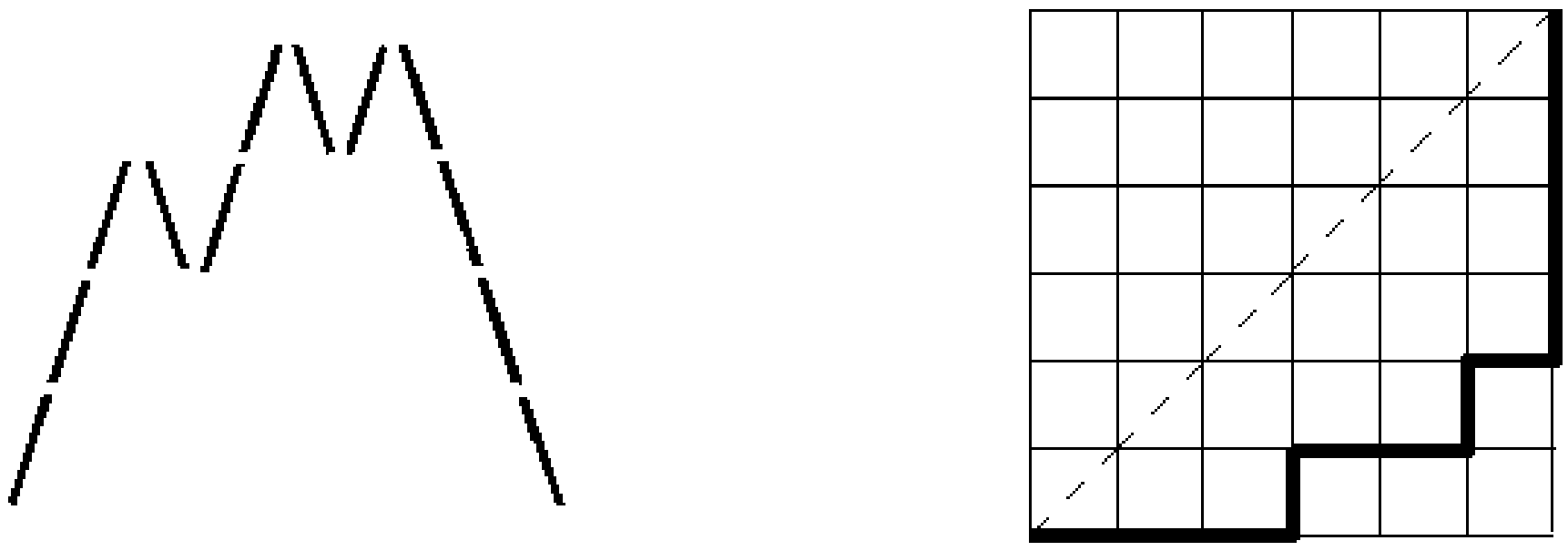

Figure 1: Corresponding Range and Path.

An $n$-by-$n$ grid is more interesting, here we see a real triangle (below the diagonal) with horizontal and vertical links of monotonic Dyck paths of length $2n$. Let's call this object a *Dyck $n$-triangle* (in our example we have the Dyck 6-triangle).

In Figure 1, both paths are identical, since they correspond to the same balanced parentheses with six open (left) and six closed (right) parentheses

(1) $$( ( ( ) ( ( ) ( ) ) ) ) .$$

We liked this definition of balanced parentheses: “a string of parentheses is valid if there are an equal number of open and closed parentheses and if you begin at the left as you move to the right, add 1 each time you pass an open and subtract 1 each time you pass a closed parenthesis, then the sum is always non-negative” [2]. In the bracket group, the left (right) parenthesis corresponds to the upstep (downstep) in the diagonal path and the horizontal (vertical) link in the monotonic path.

The bijection between diagonal paths, monotonous paths and balanced parentheses is obvious. The number of diagonal paths with $n$ upsteps and $n$ downsteps, the number of monotonic paths in a grid of $n \times n$ squares, and the number of balanced parentheses of length $2n$ are equal to the $n$-th Catalan number.

In Figure 1, both paths are drawn on a plane in different coordinate grids (the left path contains only diagonal links, in the right picture there are no such links), each grid has its own distinct coordinates. Consequently, the Dyck path can be drawn in three or even four dimensions. The question is: how many types of Dyck paths are there? Obviously, it depends on the number of possible grids and their dimensions.

# 2 Catalan lattice

**2.1.** On the plane, we usually draw two (rarely three) coordinate axes, the remaining coordinates for multidimensional objects are represented virtually. Balanced parentheses are convenient because there is direct access to all four coordinates.

Let's return to (1) and start looking at the symbols from left to right, counting the left and right parentheses. Suppose that in some $i$-th *position* (in our case $i \le 12$) we have counted $l$ *left parentheses* and $r$ *right parentheses*. Obviously, $l = r = 6$ for $i = 12$. In the general case, for the balanced parentheses of length $2n$ we have

(2) $$i = l + r, \ \ 2n \ge i \ge l \ge r \ge 0.$$

In addition, another variable is introduced, the so-called *unbalance* $j$ (the excess of the left parentheses over the right ones):

(3) $$j = l - r, \ \ n \ge l \ge j \ge 0.$$

The equalities (2) and (3) will be called a *coordinate tie*. The variables $i$ and $l$ are dominant: if $i = 0$ or $j = 0$, the others are reset. A smaller status has $j$ and $r$; $j$ and $r$ do not depend on each other. The sum or difference of $i$ and $j$ is even or zero.

The variables $i$, $j$, $l$, $r$ are non-negative integers. Let's call the *Catalan lattice* an infinite set of points in 4D with non-negative integer coordinates (the first hyperoctant) that satisfy equalities (1) and (2). Thus in the general case it can be said that the Dyck paths of arbitrary length are located in the Catalan lattice.

In Figure 1, we show the diagonal paths in the $i \times j$ grid and the monotone paths in the $l \times r$ grid. There are other versions. For example, the reader can obtain diagonal-monotonic paths in the $l \times j$ grid (diagonal upsteps and vertical downsteps). Obviously, 2D versions are modifications of the 4D Dyck triangle. We have $\binom{4}{2} = 6$ 2D modifications and $\binom{4}{3} = 4$ 3D modifications. Thus, 11 objects (including a 4D structure) are associated with the Catalan lattice. We will try to get the 4D Dyck triangle step-by-step.

**2.2.** We noted that the Dyck paths in the $l \times r$ grid are simple and convenient. On the basis of this modification, we will build Dyck triangles of greater dimension.

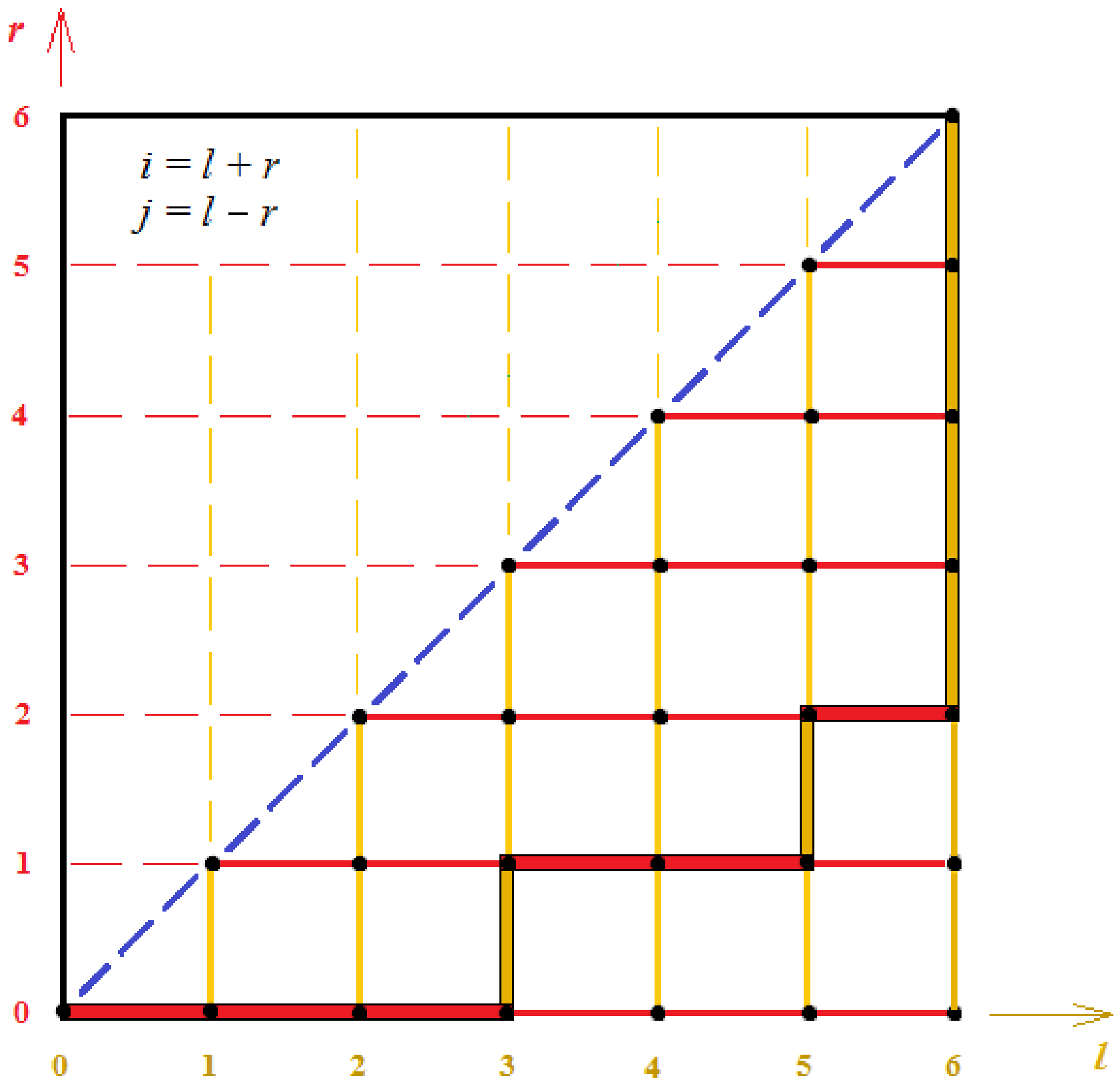

Figure 2: The Dyck triangle in a grid of 6 × 6 squares.

Let's consider in detail the monotonic path from the previous figure. In Figure 2, yellow shading shows everything related to $l$, including $l$-axis and vertical lines, *l-isolines*, which pass through points with the same coordinate $l$. Everything related to the variable $r$ is marked in red, including $r$-axis and horizontal *r-isolines*, that connect nodes with the same coordinate $r$. In a blue dotted line, we noted a diagonal line (central ray) that joins nodes for which $j = 0$ ($l = r$).

The abscissa axis is the $r$-isoline #0 (nodes with zero ordinate). Accordingly, the ordinate axis is the $l$-isoline #0 with a single accessible node (0, 0). The central ray is the $j$-isoline #0. Thus, the Dyck 6-triangle is bounded by the $r$-isoline #0 (the abscissa axis), the $l$-isoline #6 (vertical edge), and the $j$-isoline #0 (diagonal).

In Figure 2, the vertices of the Dyck triangle are follow points: the origin and the nodes (6, 0), (6, 6). Let's track the sides and vertices of the Dyck 6-triangle, going over to more complex modifications of the Catalan lattice.

**2.3.** In the $n \times n$ grid, $(n+1)(n+2)/2$ nodes is achievable. We get the same number of nodes in all modifications of the 4D Dyck $n$-triangle. Each achievable node is determined by a pair of coordinates according to the equations (2) and (3). For example, if $i$ and $j$ are given, then we get the coordinate grid $\frac{1}{2}(i+j) \times \frac{1}{2}(i-j)$.

Earlier we talked about the bijection between different types of Dyck paths of fixed length. In this case we have a *bijection* between nodes of the Dyck $n$-triangle in different modifications. It is logical to expect that the 4D triangle and its 3D modifications are as flat as the 2D projections. We will check this further.

## 3 Grid of $n \times n \times n$ cubes

In Figure 2 the square is converted to a 3D object if we add the axis $i$ or $j$. The range of values for $j$, $l$, and $r$ is the same, so in the $j \times l \times r$ grid we get the usual more convenient cube. Figure 3 shows the transition to the $j \times l \times r$ grid.

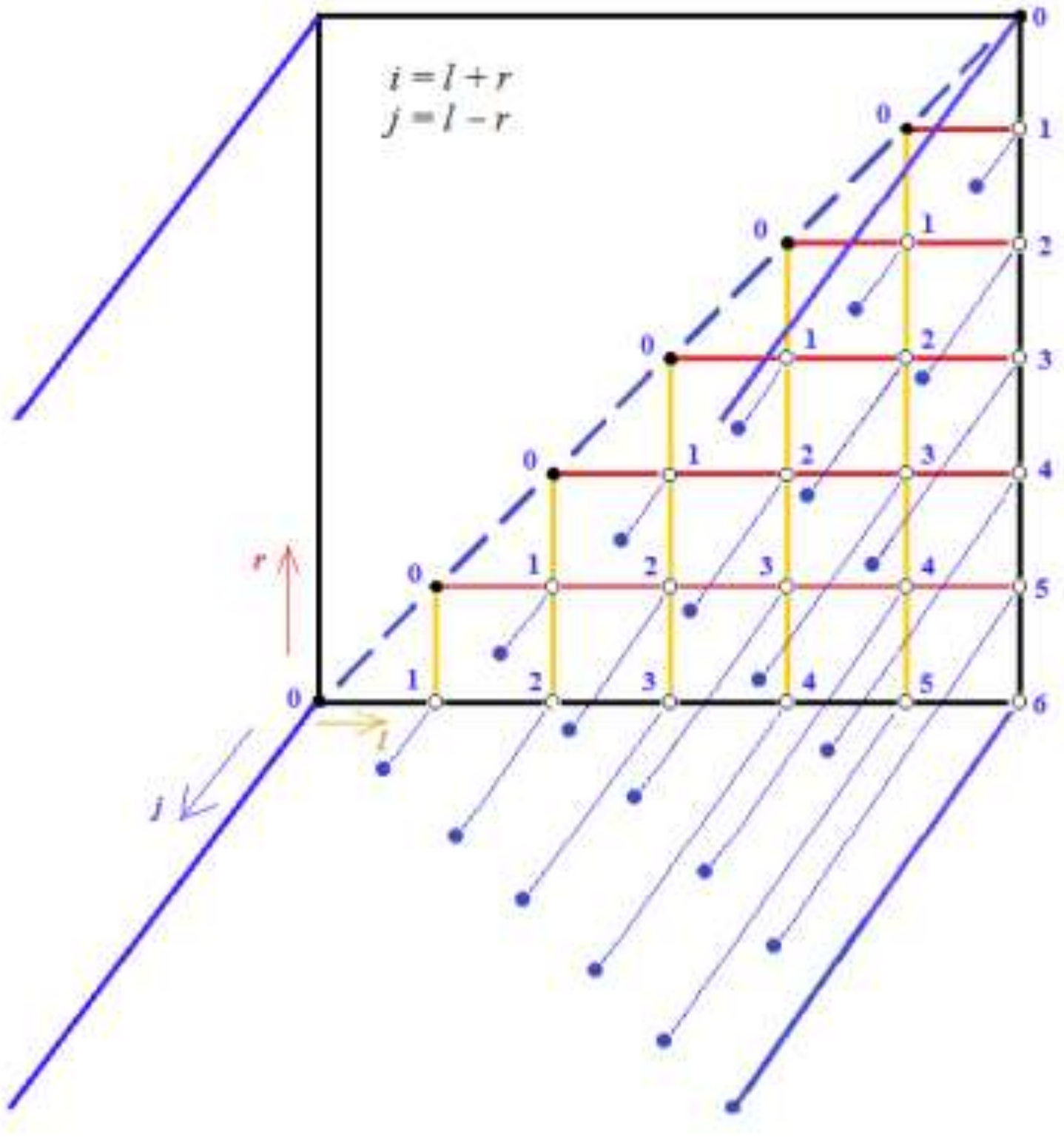

Figure 3: The transition from a 6 × 6 square to a 6 × 6 × 6 cube.

Everything related to the *j*-axis is shown in blue. The new edges of the formed cube are parallel to the *j*-axis and are drawn from the vertices of the original square.

The achievable nodes of the original square are labeled by the value of *j* according to equality (3). The blue color also shows the movement of the square nodes inside the cube, the depth of movement corresponds to the value of the labels. Zero marks the diagonal nodes of the square (zero *j*-isoline), and these nodes do not move inside the cube under construction.

It remains to draw the front face of the cube (parallel to the original square) and also *l*-isolines and *r*-isolines inside the cube. In addition, we repeat the previous Dyck path in the cube. The result is shown in Figure 4.

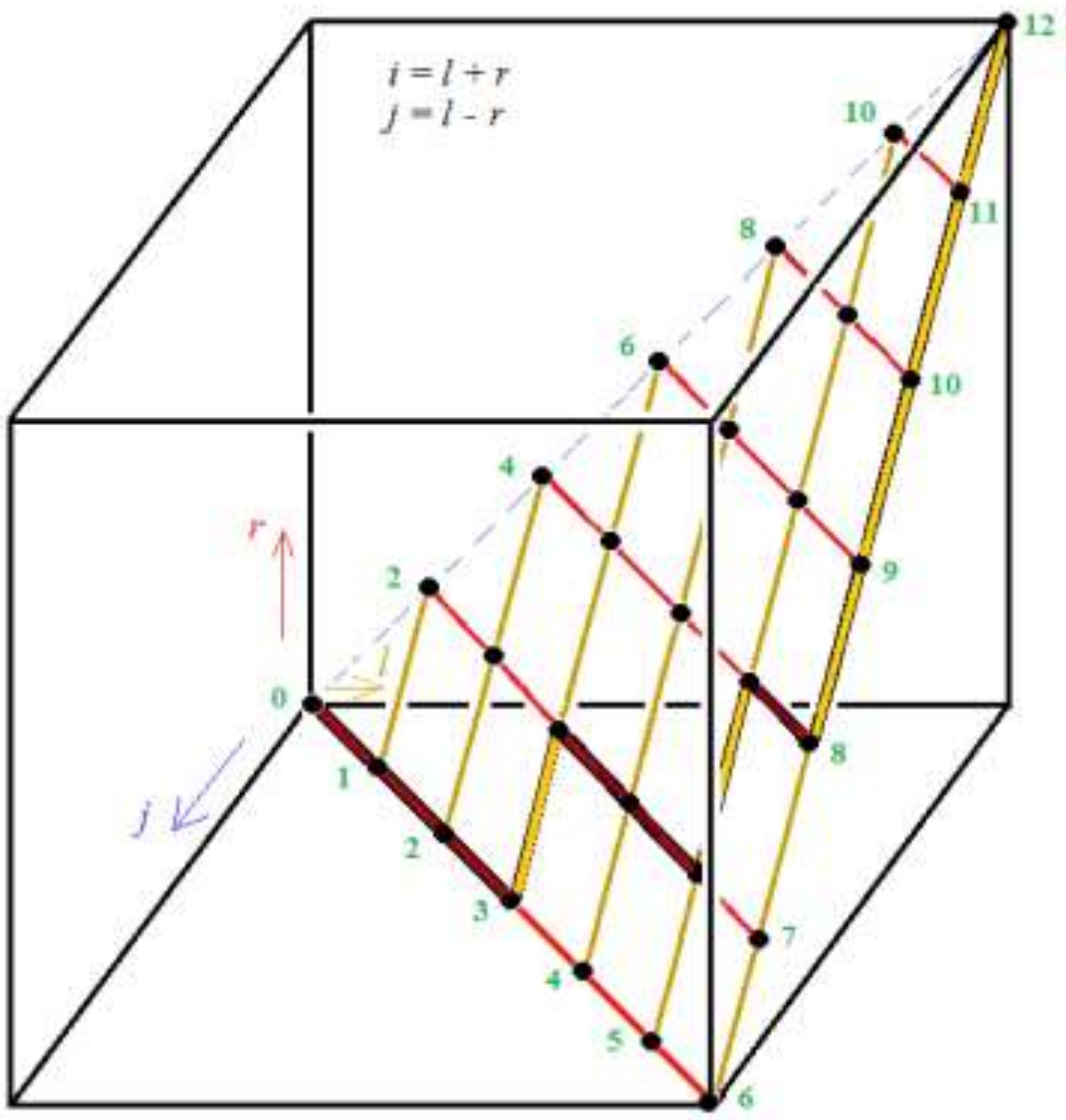

Figure 4: The Dyck 6-triangle in a grid of 6 × 6 × 6 cubes.

When you pass from a square to a cube, each square edge is converted to a cube face. In this connection, the horizontal and vertical sides of the Dyck triangle become diagonals of the cube faces. In the first octant of the $j \times l \times r$ grid, two diagonals of faces are drawn from the origin (0, 0, 0): a red ray (zero *r*-isoline) on the plane *jl* and a blue dotted ray (zero *j*-isoline) on the plane *lr*. All points between these two rays form the Catalan lattice in the 3D projection.

The Dyck *n*-triangle is flat inside an $n \times n \times n$ cube. In the Catalan lattice, this triangle is cut off by a segment of the *l*-isoline #*n* in the face that is parallel to the plane *jr*. The extreme nodes of the segment are $(n, n, 0)$ and $(0, n, n)$.

Figure 4 is not very burdened with information, so we additionally marked external nodes of the Dyck 6-triangle with the value of $i$ (the sum of $l + r$). Everything connected with the coordinate $i$ is shown in green.

# 4 Dyck triangle in double tesseract $2n \times n^3$

**4.1.** By adding the $i$-axis we get a 4D space. Among the four-dimensional figures, the most popular is the tesseract (hypercube, tetracube, octagonal) which is obtained if we move the ordinary cube into the fourth dimension by an additional edge and then connect the vertices of both cubes. In our case, the fourth edge has double length, so we are dealing with a four-dimensional parallelepiped in the first hyperoctant (or orthant). Let's call it a *double tesseract*.

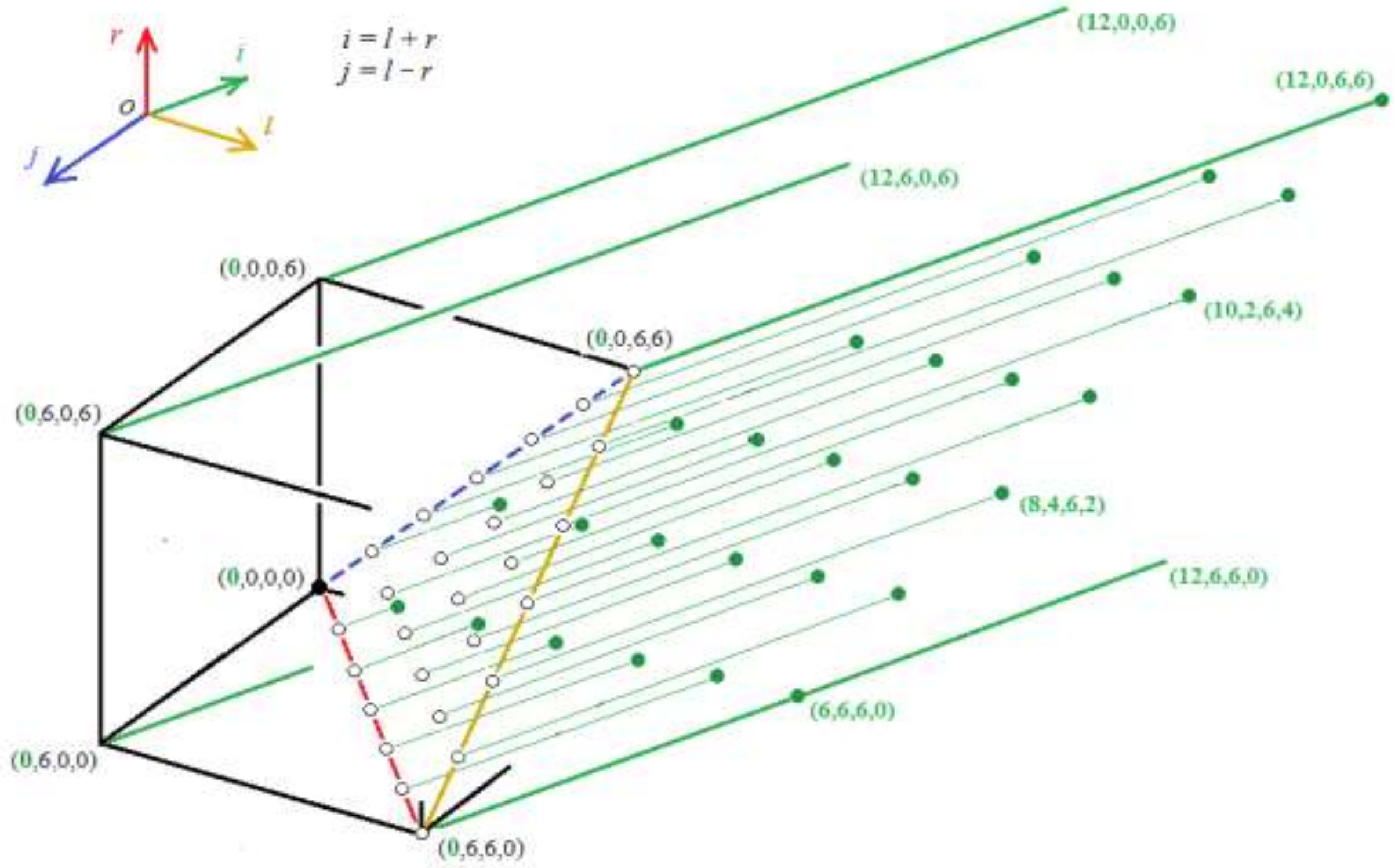

Figure 5: Transition to the $12 \times 6^3$ grid.

Figure 5 shows the transition from the 3D cube (the first octant in the $j \times l \times r$ grid) to the double tesseract (the first hyperoctant in the $i \times j \times l \times r$ grid). The direction of axes is shown at the top left. We placed the origin in the far bottom node of the source cube. Green indicates that it is associated with the 4th dimension.

For a better view, some nodes are cut. We showed the movement of achievable nodes of the original cube inside the double tesseract according to the value of $i$ (see Figure 4). Only one achievable node, the origin, does not move.

**4.2.** In Figure 6, the $i \times j \times l \times r$ grid shows the resulting double tesseract, inside the 4D Dyck triangle is located (familiar Dyck path is repeated). The triangle inside is darkened. The individual nodes of the tesseract cut again.

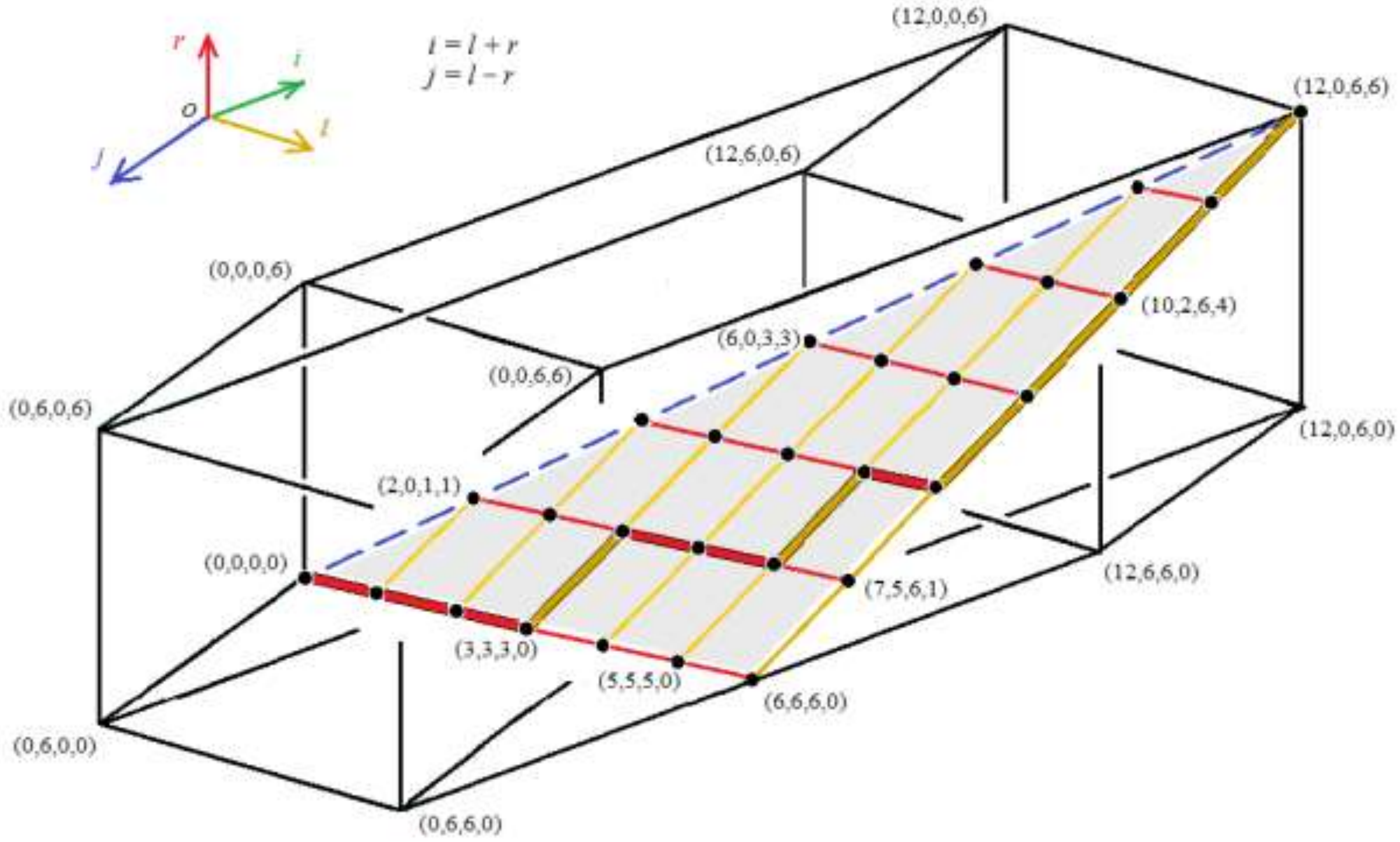

Figure 6: The Dyck 6-triangle in a double tesseract of $12 \times 6^3$ 4-cubes.

A double tesseract (like a usual tesseract) has 16 vertices, 32 edges, and 8 three-dimensional faces. Faces are formed by the intersection of a tesseract with eight hyperplanes. In a normal tesseract, all 3D faces are cubes. In our case, only two faces are 3D cubes, and this is clearly seen in Figure 6.

**4.3.** In Figure 6, the four-dimensional Dyck 6-triangle and all its Dyck paths begin at the origin and end at the node (12, 0, 6, 6). Both these nodes and five more accessible nodes form the *blue side* of the Dyck 6-triangle ($j$-isoline #0, blue dotted line). In addition, we marked the *red side* ($r$-isolines #0) and the *yellow side* ($l$-isoline #6) of the triangle. All three sides are straight lines in the 4D space.

In the general case, the blue side of the Dyck $n$-triangle is the diagonal of the regular parallelepiped (the 3D face of the double tesseract $2n \times n^3$), which contains all nodes with $j = 0$. The parallelepiped is bound by the following nodes:

$(0,0,0,0)$, $(0,0,0,n)$, $(0,0,n,n)$, $(0,0,n,0)$, $(2n,0,n,0)$, $(2n,0,0,0)$, $(2n,0,0,n)$, $(2n,0,n,n)$.

In this list, the origin and the last point are the extreme nodes of the blue side.

The other two sides of the Dyck $n$-triangle intersect at the point $(n, n, n, 0)$, the center of the edge of the double tesseract with the extreme points $(0, n, n, 0)$ and ($2n$, $n$, $n$, 0). The red side is placed in the regular parallelepiped (another 3D face), which contains all points with $r = 0$. The parallelepiped is bounded by points

$(0,0,0,0), (0,n,0,0), (0,n,n,0), (0,0,n,0), (2n,0,n,0), (2n,0,0,0), (2n,n,0,0), (2n,n,n,0).$

Let us divide this parallelepiped into two equal 3D cubes. Then in the left cube, the red side of the $n$-triangle becomes a diagonal in the form of a ray with points $(x, x, x, 0)$, $x \le n$ .

The yellow side of the Dyck $n$-triangle is placed in a parallelepiped (another 3D face of the double tesseract), which contains all points of a double tesseract with $l = n$. This parallelepiped is bounded by points

$(0,0,n,0), (0,n,n,0), (0,n,n,n), (0,0,n,n), (2n,0,n,n), (2n,0,n,0), (2n,n,n,0), (2n,n,n,n).$

Divide this parallelepiped into two equal cubes, then in the right cube, the yellow side of Dyck triangle will become a diagonal.

The straightness of the sides of the triangle can be verified by comparing the lengths of the vectors in the 4D space. For example, in the blue side there are six links of the same length $\sqrt{(2^2 + 0^2 + 1^2 + 1^2)} = \sqrt{6}$. On the other hand, the entire length of the blue side is determined by the coordinates of the point (12, 0, 6, 6) and is equal to $\sqrt{(12^2 + 0^2 + 6^2 + 6^2)} = 6\sqrt{6}$.

Also, the red side includes six links with a length of $\sqrt{(1^2 + 1^2 + 1^2 + 0^2)} = \sqrt{3}$. The length of the entire red side is determined by the coordinates of the last vertex (6, 6, 6, 0) and is equal to $\sqrt{(6^2 + 6^2 + 6^2 + 0^2)} = 6\sqrt{3}$. We obtain the same lengths for the yellow side, it follows that the 4D Dyck 6-triangle is not only flat, but also isosceles.

In Figure 2, the red side and yellow side are perpendicular. Let us verify this in 4D space. In the Dyck $n$-triangle, at the intersection of the red side and the yellow side (see Figure 6), we select three points $A(n-1, n-1, n-1, 0)$, $B(n, n, n, 0)$, and $C(n+1, n-1, n, 1)$. For line AB and line BC, the direction vectors are (1, 1, 1, 0) and (1, –1, 0, 1) respectively. It remains to calculate the scalar product of the directing vectors: $1\times1 + 1\times(-1) + 1\times0 + 0\times1 = 0$. This proves the perpendicularity of AB and BC and, respectively, the perpendicularity of the red side and yellow side of the Dyck $n$-triangle.

**Conclusion.** Thus, the four-dimensional Dyck $n$-triangle is a flat right-angled isosceles triangle. It is interesting that a tie of coordinates ensures the independence of modifications. Each modification is self-contained and is able to reincarnate into other modifications (changing coordinate grids) up to the 4D triangle.

We obtained a three-dimensional picture of a double tesseract by moving the ordinary cube to the fourth dimension. But there is another popular projection of a tesseract into 3D space. Figure 7 shows the well-known Schlegel diagram, which is two nested cubes; the neighboring vertices of cubes are connected by straight lines.

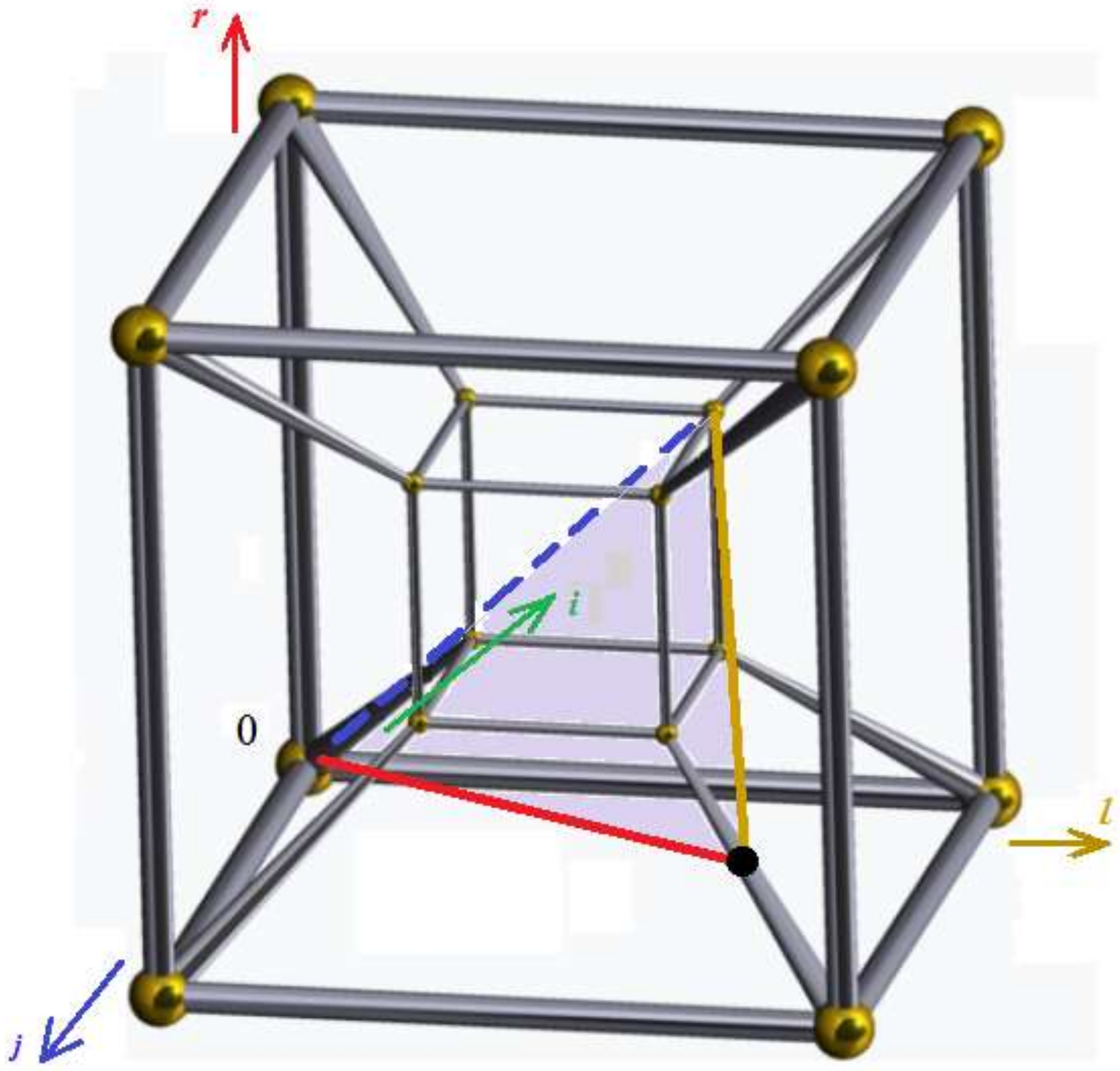

Figure 7: The Dyck *n*-triangle in a Schlegel diagram.

In a 4D space, the outer and inner cubes have the same dimensions. In our case, both cubes are connected by edges of twice the length, and the result is a double tesseract. We placed the origin in the lower node of the outer cube, and the three axes *j*, *l*, and *r* look logical. The fourth virtual *i*-axis "goes" inside both cubes.

The 4D Dyck *n*-triangle is shown in the diagram by the following three points: (1) the origin, (2) the far right point of the top face of the inner cube, and (3) the center of the doubled edge opposite the origin. The red and blue sides of the Dyck triangle are located on the rays, which limit the Catalan lattice in the 4D space.

**Acknowledgements.** I wish to thank Igor Pak (Mathematics Department at UCLA, USA) for the detailed and comprehensive historical review of Catalan-like numbers, and thanks for a selection of recent works on the subject. Special thanks to Bruce Sagan (Michigan State University, USA) who always pays great attention to my work.

Gzhel State University, Moscow, 140155, Russia

http://www.en.art-gzhel.ru/